\documentclass[10pt,twoside]{article}
\usepackage{graphicx}
\usepackage{amsmath,amsfonts}
\usepackage{Latex-document}
\usepackage{epsf}

\title{\bf Singular Combinatorics\thanks{Work supported in part by the IST Programme of the EU under
contract number IST-1999-14186 (ALCOM-FT).}\vskip 6mm}
\author{Philippe Flajolet\vspace*{-0.5cm}\thanks{Algorithms Project,
INRIA-Rocquencourt, 78153 Le Chesnay, France. E-mail:
Philippe.Flajolet@inria.fr}}
\date{\vspace{-8mm}}

\newtheorem{theorem}{Theorem}
\def\implies{\Longrightarrow}
\def\Z{\mathbb Z}
\def\cal{\mathcal}
\def\frak{\mathfrak}
\def\hat{\widehat}
\def\ds{\displaystyle}
\def\GF{{\sc gf}}
\def\OGF{{\sc ogf}}
\def\EGF{{\sc egf}}
\def\Arg{\operatorname{Arg}}
\def\Li{\operatorname{Li}}

\newcounter{exno}

\newenvironment{example}%
        {\refstepcounter{exno}%
        \medbreak\noindent%
        \setlength{\baselineskip}{0.85\baselineskip}%
        \begin{small}%
    \setlength{\abovedisplayskip}{3.0ptplus0ptminus0.5pt}%
    \setlength{\belowdisplayskip}{3.0ptplus0ptminus0.5pt}%
    \setlength{\abovedisplayshortskip}{1.5ptplus0ptminus0pt}%
    \setlength{\belowdisplayshortskip}{1.5ptplus0ptminus0pt}%
    \predisplaypenalty=-100\postdisplaypenalty=0%
        {\bf Example}~{\hbox{\bf\arabic{exno}.}}}%
        {\end{small}\hfill{\normalsize $\Box$}\medbreak}

\newenvironment{proof}{{\bf Proof.}}{\hfill$\Box$\par\smallskip}
\begin{document}

\maketitle \thispagestyle{first} \setcounter{page}{561}

\begin{abstract}\vskip 3mm
Combinatorial enumeration leads to counting generating
functions presenting a wide variety of analytic types. Properties of
generating functions at singularities encode valuable information regarding
asymptotic counting and limit probability distributions
present in large random structures. ``Singularity analysis'' reviewed here
provides constructive estimates that are applicable
in several areas of combinatorics. It constitutes a
complex-analytic Tauberian procedure by which combinatorial
constructions and asymptotic--probabilistic laws can be
systematically related.

\vskip 4.5mm


\end{abstract}

\vskip 12mm

\section{Introduction}

\vskip-5mm \hspace{5mm}

Large random combinatorial structures tend to exhibit great statistical
regularity. For instance, an overwhelming proportion
of the graphs of a given large size are connected,
and a fixed pattern
is almost surely contained in a long random string, with its number
of occurrences satisfying central and local limit laws.
The objects considered (typically, words, trees, graphs,
or permutations) are given by
construction rules of the kind classically studied by combinatorial
analysts via \emph{generating functions} (abbreviated as {\GF}s).
A fundamental problem is then to extract
asymptotic information on coefficients of a {\GF}
either explicitly given by a formula or implicitly determined by a
functional equation. In the univariate case, asymptotic counting
estimates are derived; in the multivariate case, moments and
limit probability laws of characteristic parameters
will be obtained.

In what follows, given a combinatorial class $\mathcal C$,
we let $C_n$ denote the number of objects in~$\mathcal C$ of size~$n$ and
introduce the \emph{ordinary} and \emph{exponential} {\GF}
({\OGF}, {\EGF}),
\[
\hbox{{\OGF}:}\quad
C(z):=\sum_{n\ge0} C_n z^n,\qquad
\hbox{{\EGF}:}\quad
\widehat{C}(z):=\sum_{n\ge0} C_n \frac{z^n}{n!}.
\]
Generally, {\EGF}s and {\OGF}s serve for the enumeration of labelled classes
(atoms composing objects are distinguished by
labels) and unlabelled classes, respectively.
One writes $C_n=[z^n]C(z)=n!\,[z^n]\,\widehat C(z)$,
with~$[z^n](\cdot)$ the coefficient extractor.

General rules for deriving
{\GF}s from combinatorial specifications
have been widely developed by various schools
starting from the 1970's and these lie at the heart of contemporary
combinatorial analysis.
They are excellently surveyed in books of
Foata \& Sch\"utzenberger (1970), Comtet (1974),
Goulden \& Jackson
(1983), Stanley (1986, 1998),  Bergeron, Labelle \& Leroux (1998).
We shall retain here the following simplified scheme
relating combinatorial constructions and operations over {\GF}s:
\begin{equation}\label{basicdic}
\hbox{\small
\begin{tabular}{ll|ll}
\hline
\em Construction & & \em Labelled case & \em Unlabelled case  \\
\hline\hline
Disjoint union & $\cal F +\cal G$ &  $f(z)+g(z)$ & $\hat f(z)+\hat g(z)$  \\
Product & $\cal F\times\cal G$, $\cal F\star\cal G$ & $f(z)\cdot g(z)$
& $\hat f(z)\cdot \hat g(z)$  \\
Sequence & $\frak S\{\cal F\} $ & $(1-f(z)))^{-1}$ &$(1-f(z)))^{-1}$ \\
Set & $\frak P\{\cal F\}$ & $\exp(f(z)$)& $\exp\left(f(z)+\frac12f(z^2)+\cdots\right)$  \\
Cycle &  $\frak C\{\cal F\}$ & $\log(1-f(z))^{-1}$ & $\log(1-f(z))^{-1}+
\cdots$ \\
\hline
\end{tabular}
}
\end{equation}

Such operations on {\GF}s
yield a wide variety of analytic functions, either  given explicitly or
as solutions to functional equations in the case of recursively defined
classes. It is precisely the goal of singularity analysis to
provide means for extracting asymptotic informations.
What we termed ``singular combinatorics''
aims at relating combinatorial form
and asymptotic-probabilistic form by exploiting
\emph{complex-analytic} properties
of generating functions.
Classical approaches~\cite{Bender74}
are Tauberian theory and Darboux's method,
an offspring of elementary Fourier analysis largely developed
by P\'olya for his programme of combinatorial chemistry~\cite{Polya37}.
The path followed here, called ``singularity analysis'' after~\cite{FlOd90b},
consists in developing a systematic correspondence between the local
behaviour of a function near its singularities and the asymptotic form
of its coefficients. (An excellent survey of central aspects of the theory
is offered by Odlyzko in~\cite{Odlyzko95}.)

\section{Basic singularity analysis}\setzero
\vskip-5mm \hspace{5mm}

Perhaps the simplest coefficient estimate is
$[z^n](1-z)^{-\alpha}\sim {n^{\alpha-1}}/{\Gamma(\alpha)}$,
a consequence of the binomial expansion and Stirling's formula.
For the basic scale 
\[
\sigma_{\alpha,\beta}(z) = (1-z)^{-\alpha}\left(\frac1z \log(1-z)^{-1}
\right)^\beta,
\]
much more is available and one has a fundamental translation
mechanism~\cite{FlOd90b}:

\begin{figure}\small
\vspace*{2.5truemm}
\begin{center}\setlength{\unitlength}{0.6truemm}
\begin{picture}(140,40)(15,0)
\put(0,16){
\begin{picture}(60,16)(-8,-8)
\put(0.5,-4){\bf 1}
\put(0,0){\line(1,0){40}}
\thicklines
\put(40,-8){\vector(-1,0){40}}
\put(0,8){\vector(1,0){40}}
\put(0,0){\oval(16,16)[l]}
\thinlines
\put(35,0){\vector(0,1){8}}
\put(35,8){\vector(0,-1){8}}
\put(40,3){$\frac1n$}
\put(20,10){$\mathcal H$}
\end{picture}}
\put(75,0){
\begin{picture}(70,40)(8,-24)
\put(0.5,-4){\bf 1}
\put(0,0){\line(1,0){40}}
\thicklines
\put(40,-21.33){\vector(-3,1){40}}
\put(0,8){\vector(3,1){40}}
\put(0,0){\oval(16,16)[l]}
\thinlines
\put(0,0){\vector(0,1){8}}
\put(0,8){\vector(0,-1){8}}
\put(0,3){\ $\frac1n$}
\put(30,8){\vector(-1,2){4}}
\put(30,8){\vector(0,-1){8}}
\put(25,7){$\vartheta$}
\put(6,16){$\mathcal H(\vartheta)$}
\end{picture}}
\put(110,-4){{\ \epsfxsize4truecm\ \epsfbox{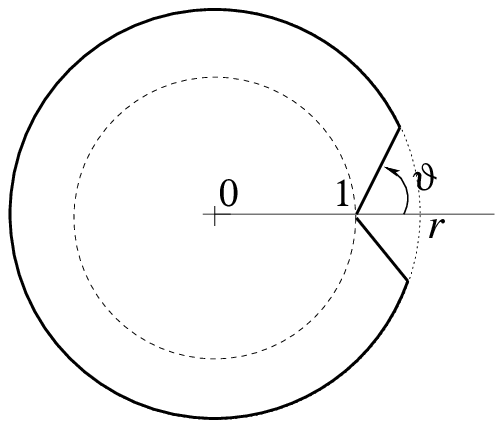}}}
\end{picture}
\end{center}
\vspace*{-0.72truecm}
\caption{\label{hankel-figs}The Hankel contours, $\mathcal H$ and
$\mathcal H(\vartheta)$, and a $\Delta$-domain.}
\end{figure}

\begin{theorem}[Coefficients of the basic scale]\label{thm1}
For $\alpha\in\mathbb C\setminus\mathbb Z_{\le 0}$,
$\beta\in\mathbb C$, one has
\begin{equation}\label{sigmabn}
[z^n]\, \sigma_{\alpha,\beta}(z)\mathop{\sim}_{n\to\infty}
\frac{n^{\alpha-1}}{\Gamma(\alpha)}\left(\log n\right)^\beta.
\end{equation}
\end{theorem}
\begin{proof}
The estimate is derived starting from Cauchy's coefficient
formula,
\[
[z^n]f(z)\sim \frac{1}{2i\pi}\int_{\mathcal\gamma}
f(z)\, \frac{dz}{z^{n+1}},
\]
instantiated with $f=\sigma_{\alpha,\beta}$. The idea is then
to select for~$\gamma$
a contour $\mathcal{H}$ that is of Hankel type and follows the half-line
$(1,+\infty)$ at distance exactly $1/n$ (Fig.~\ref{hankel-figs}). This
superficially resembles
a saddle-point contour, but with the integral normalizing to
Hankel's representation of the Gamma function, hence the factor
$\Gamma(\alpha)^{-1}$.
\end{proof}
The method of proof is very flexible: it applies for instance
to  iterated logarithmic terms ($\log\log$)
while providing full
asymptotic expansions; see~\cite{FlOd90b} for details.

A remarkable fact, illustrated by Theorem~\ref{thm1}, is that larger functions
near the singularity $z=1$ give rise to larger coefficients
as $n\to\infty$. This is a general phenomenon
under some suitable auxiliary conditions, expressed here in
terms of \emph{analytic continuation}: a $\Delta$-domain is an
indented disc defined by ($r>1$, $\vartheta<\pi/2$)
\[
\Delta(\vartheta,r):=\left\{ \, z~\bigm|~|z|<r,~
\vartheta<\Arg(z-1)<2\pi-\vartheta,~z\not=1\,\right\}.
\]

\begin{theorem} [{\boldmath $O$}-transfer] \label{thm2}
With $f(z)$ continuable to a $\Delta$-domain and
$\alpha\not\in \Z_{\le0}$:
\[
f(z)\mathop{=}_{z\to1,~ z\in\Delta} O\left(\sigma_{\alpha,\beta}(z)\right)
\qquad
\implies\qquad
[z^n]\, f(z) =O\left( [z^n]\, \sigma_{\alpha,\beta}(z)\right).
\]
\end{theorem}
\begin{proof}
In Cauchy's coefficient formula, adopt an integration contour
$\cal H(\vartheta)$
passing at distance $1/n$ left of the singularity~$z=1$, then
escaping outside of the unit disc within~$\Delta$. Upon setting $z=1+t/n$,
careful approximations yield the result~\cite{FlOd90b}.
\end{proof}

This theorem allows one to transfer
error terms in the asymptotic
expansion of a function at its singularity (here $z=1$)
to asymptotics of the coefficients.
The Hankel contour technique is quite versatile
and a statement similar to Theorem~\ref{thm2} holds with
$o(\cdot)$-conditions.
replacing $O(\cdot)$-conditions.
The case of~$\alpha$ being a negative integer is
covered by minor adjustments due to $1/\Gamma(\alpha)=0$;
see~\cite{FlOd90b}.
In concrete terms: \emph{Hankel contours combined with Cauchy
coefficient integrals accurately ``capture''
the singular behaviour of a function.}

By Theorems~\ref{thm1} and~\ref{thm2},
whenever a function admits an asymptotic expansion near $z=1$
in the basic scale,
one has the implication, with $\sigma \succ \tau\succ \ldots\succ\omega$,
\[
f(z)=\lambda \sigma(z)+\mu\tau(z)+\cdots +O(\omega(z))
\quad
\implies
\quad
f_n=\lambda\sigma_n+\mu\tau_n+\cdots+O(\omega_n),
\]
where $f_n=[z^n]f(z)$. In other words, a
\emph{dictionary}
translates singular expansions of functions into
the asymptotic forms of coefficients. Analytic continuation
and validity
of functions' expansions outside of the unit circle is
a logical necessity, but once granted,
application of the method becomes quite mechanical.

For combinatorics, singularities need not be placed at~$z=1$. But since
$[z^n]f(z)\equiv \rho^{-n}[z^n]f(z\rho)$,
the dictionary can be used for (dominant) singularities
that lie anywhere in the complex plane. The case of finitely many
dominant singularities can also be dealt with
(via composite Hankel contours) to
the effect that the translations of local singular expansions get composed
additively. In summary one has from function to coefficients:

\medskip

\centerline{\small \fbox{location+nature of singularity (fn.) $\implies$ exponential+polynomial
asymptotics (coeff.)}}

\smallskip

\begin{example} \emph{2-Regular graphs} (Comtet, 1974). The class $\cal G$
of (labelled) 2-regular graphs can be specified as
sets of unordered cycles each of length at least~$3$. Symbolically:
$$
\cal G \cong \frak P \{ \frac12\frak C_{\ge3} \{ \cal Z \}\}
\quad\hbox{so that}\quad
\hat G(z)=\exp\left(\frac12\log(1-z)^{-1}-\frac{z}{2}
-\frac{z^2}{4}\right)
=\frac{e^{-z/2-z^2/4}}{\sqrt{1-z}}
.$$
$\cal Z$ represents a single atomic
node. The function $\hat G(z)$ is
singular at~$z=1$, and
$$
\cal G(z)
\mathop{\sim}_{z\to1} e^{-3/4}(1-z)^{-1/2}
\qquad\implies \qquad
\frac{G_n}{n!}\mathop{\sim}_{n\to\infty} \frac{e^{-3/4}}{\sqrt{\pi n}}.
$$
This example can be alternatively treated
by Darboux's method. 
\end{example}

\begin{example}
The \emph{diversity index} of a tree (Flajolet, Sipala \& Steyaert, 1990)
is the number of non-isomorphic terminal subtrees, a quantity
also equal to the size of maximally compact representation of the tree
as a directed  acyclic graph and related to common subexpression
sharing in computer science applications.
The mean index of a random binary tree of size~$2n+1$
is asymptotic to $Cn/\sqrt{\log n}$,
where $C=\sqrt{8\log 2/\pi}$.  This results from
an exact {\GF} obtained by inclusion-exclusion:
\[
K(z)=\frac{1}{2z}\sum_{k\ge0} \frac{1}{k+1}\binom{2k}{k}
\left(\sqrt{1-4z+4z^{k+1}}-\sqrt{1-4z}\right).
\]
Singularities accumulate geometrically to the right of $1/4$ while~$K(z)$ is
$\Delta$--continuable. The
unusual singularity type ($1/\sqrt{X\log X}$)
precludes the use of Darboux's method.
\end{example}

Rules like those of Table~(\ref{basicdic})
preserve analyticity and analytic continuation.
Accordingly, generating functions associated with combinatorial
objects described by simple construction rules usually have {\GF}s
amenable to singularity analysis. The method is systematic enough, so that
an implementation within computer algebra systems is even possible
as was first demonstrated by Salvy~\cite{Salvy91}.

\section{Closure properties}

\vskip-5mm \hspace{5mm}

In what follows, we say that a function is
\emph{amenable to singularity analysis}, or
``\emph{of S.A.~type}'' for short,
if it is $\Delta$-continuable and admits there a singular expansion in the
scale $\cal S=\{\sigma_{\alpha,\beta}(z)\}$.
First, functions of S.A.~type include polylogarithms:
\begin{theorem} The generalized polylogarithms~$\Li_{\alpha,k}$
are of S.A. type, where
$\Li_{\alpha,k}(z):=\sum_{n\ge1}{(\log n)^k}{n^{-\alpha}}z^n$,
$k\in\mathbb Z_{\ge0}$.
\end{theorem}

The proof makes use of the Lindel\"of representation
\[
\sum_{n\ge1}\phi(n)(-z)^{n}
=-\frac{1}{2i\pi}\int_{1/2-\infty}^{1/2+\infty}
\phi(s) z^s \frac{\pi}{\sin \pi s}\, ds,
\]
in a way already explored by Ford~\cite{Ford60sp}, with Mellin transform
techniques providing validity of the singular expansion in a
$\Delta$-domain~\cite{Flajolet99}.

\begin{example} {\em Entropy computations.}
The {\GF} of $\{\log(k!)\}$ is $(1-z)^{-1}\Li_{0,1}(z)$,
which is of S.A.~type.
The entropy of the binomial distribution,
$\pi_{n,k}=\binom{n}{k}p^k(1-p)^{n-k}$, results:
\[
H_{n}:=-\sum \pi_{n,k}\log \pi_{n,k}
\mathop{\sim}_{n\to\infty}
\frac12\log n+\frac12+\log\sqrt{2\pi p(1-p)}+\cdots\,
.\]

Such problems are of interest in information theory, where
redundancy estimates  precisely depend on higher order asymptotic
properties (Jacquet-Szpankowski, 1998). Full expansions for
functionals of the Bernoulli distribution are also obtained
systematically~\cite{Flajolet99}.
\end{example}

As it is well-known, asymptotic expansions can be integrated
while differentiation of asymptotic expansions is
permissible for functions analytic in sectors:
\begin{theorem}
Functions of S.A. type are closed under
differentiation and integration.
\end{theorem}

Finally, the Hadamard product of two series, $f\odot g$ is defined
as the termwise product:
$f(z)\odot g(z)=\sum_{n} f_n g_n z^n$,
if
$f(z)=\sum_{n} f_n z^n$, $g(z)=\sum_{n} f_n z^n$.
Hadamard (1898) proved 
that singularities get composed
multiplicatively.
Finer composition properties~\cite{FiFlKa02}
result from an adaptation of Hankel contours to Hadamard's formula
\[
f(z)\odot g(z)=\frac{1}{2i\pi}\int_\gamma f(t) g\left(\frac{w}{t}\right)\,
\frac{dt}{t}.
\]
\begin{theorem} Functions of S.A. type are closed under Hadamard product.
\end{theorem}
\begin{example} \emph{Divide-and-conquer} algorithms solve recursively
a problem of size~$n$ by splitting it into two subproblems and
recombining the partial solutions. Under the assumption
of randomness preservation, the expected costs~$f_n$
satisfy a ``tree recurrence'' of the form
\[
f_n=t_n+\sum_{k} \pi_{n,k}\left(f_k+f_{n-k-a}\right),\qquad a\in\{0,1\},
\]
where the ``toll'' sequence~$t_n$ usually has a simple form
(e.g., $n^\beta, \log n$) and the $\pi_{n,k}$
characterize the stochastic splitting process.
The corresponding {\GF}s then satisfy an equation
$
f(z)=t(z)+\cal L \left[ f\right](z)
$,
where the linear  operator $\cal L$ reflects the splitting probabilities.
For instance, binary search trees and the Quicksort algorithm
have $\cal L[f(z)]=2\int_0^z f(x) dx/(1-x)$. One then has
in operator notation
$
f(z)= (I-\cal L )^{-1}[t](z),
$
where the quasi-inverse acts as a ``singularity transformer''.
Closure theorems allow for an asymptotic classification of the cost functions
induced by various tolls under various probabilistic models
mirrored  by the splitting probabilities~\cite{FiFlKa02}.
\end{example}

\section{Functional equations}

\vskip-5mm \hspace{5mm}

\emph{Algebraic functions} have expansions at singularities that are
expressed by fractional power series (Newton-Puiseux). Consequently,
they are of S.A.~type with rational exponents; accordingly
their coefficient expansions are linear
combinations of algebraic elements of the form~$\omega^n n^{r/s}$,
with $\omega$ and algebraic number and $r/s\in\mathbb Q$.
By Weierstrass preparation, such properties extend to
many \emph{implicit} {\GF}s. For instance, {\GF}s of combinatorial families of
trees constrained to have degrees in a finite set have a
branch point of type $\sqrt{X}$ at their dominant singularity,
which in turn corresponds to the ``universal'' asymptotic form $T_n\sim c
\omega^n n^{-3/2}$ for coefficients (P\'olya, 1937; Otter 1948; Meir-Moon, 1978).

Order constraints in labelled classes are known to correspond to
integral operators and, in the recursive case, there result  {\GF}s
determined by
\emph{ordinary differential equations}.
Furthermore, moments (i.e., cumulative values) of
additive functionals of combinatorial structures
defined by recursion and order constraints have {\GF}s that
satisfy \emph{linear
differential equations}, for which there is a well-established
classification theory
going back to the nineteenth century. In particular, in
the \emph{Fuchsian case}, singularity analysis
applies unconditionally, so that the resulting coefficient estimates are
linear combinations of terms $\omega^n n^\alpha (\log n)^k$, with
$\omega,\alpha$ algebraic and $k\in\mathbb Z_{\ge0}$. This covers
a large subset of the class of ``holonomic'' functions, of
which Zeilberger has extensively demonstrated the
expressive power in combinatorial analysis~\cite{WiZe85b,Zeilberger90}.

\begin{example}\label{quad-exa}  {\em Quadtrees} are a way to superimpose a
hierarchical partitioning on sequences of
points in~$d$-dimensional space: the first point is taken as the root
of the tree and it partitions the whole space in~$2^d$ orthants
in which successive points are placed and then made to refine the
partition~\cite{Mahmoud92}. The problem is expressed by a linear
differential equation
with coefficients in $\mathbb{C}(z)$.
The average cost of finding a point knowing
only $k$ of its coordinates
is of the asymptotic form $c\, n^\alpha$ with~$\alpha=\alpha(k,d)$ an
algebraic number of degree~$d$. For instance $k=1$ and~$d=2$ yield
a solution involving a ${}_2F_1$ hypergeometric function
as well as $\alpha=
(\sqrt{17}-3)/2\doteq 0.56155$, in contrast to an
exponent~$\frac12$ that would correspond to a perfect partitioning,
i.e., a regular grid (Flajolet, Gonnet, Puech \& Robson, 1993).
\end{example}



\emph{Substitution equations} correspond to ``balanced structures''
of combinatorics.
An important r\^ole in the development of the theory has been played
by Odlyzko's analysis~\cite{Odlyzko82} of 2-3 trees (such trees have internal nodes of
degree 2 and 3 only and leaves are all at the same depth). The
{\OGF} satisfies the equation
$
T(z)=z+T(\tau(z))$, with $\tau(z):=z^2+z^3,
$
and has a singularity at~$z=1/\phi$, a fixed point of~$\tau$,
($\phi=(1+\sqrt{5})/2$).
The singular expansion involves periodic oscillations, corresponding
to infinitely many singular exponents having a common real part.
Singularity analysis extends to this case and the number of balanced
2-3 trees is found to be of the form
$
T_n\sim \frac1n \phi^n \Omega(\log n),
$
for some nonconstant smooth periodic function~$\Omega$.

A similar problem of \emph{singular iteration}
arises in the analysis of the height of binary
trees~\cite{FlOd82}. The {\GF}~$y_h(z)$ of trees of height at most~$h$
satisfies the Mandelbrot recurrence
$y_h = z+y_{h-1}^2$, with
$y_0=z$.
The fixed point is the {\GF} of binary trees, that is, of Catalan numbers,
$y_{\infty}=(1-\sqrt{1-4z})/2$ which has its dominant singularity
at~$1/4$. The analysis of moments of the distribution of height
turns out to be equivalent to developing \emph{uniform} approximations to
$y_h(z)$ as $z\to1/4$ and $h\to\infty$ simultaneously, this
for~$z$ in a
$\Delta$-domain. The end result, by singularity analysis and the
moment method, is: {\em the height of a random binary tree
with~$n$ external nodes when normalized by a factor of $1/(2\sqrt{n})$
converges in distribution to a theta law defined by the density}
$
4x\sum_{k\ge1} k^2(2k^2x2-3)e^{-k^2x^2}.
$
The result extends to all simple families of trees in the sense of
Meir and Moon and it provides pointwise estimates of the
proportion of trees of given height, that is, a \emph{local limit law}.

Generalized digital trees (Flajolet \& Richmond, 1992)
correspond to a \emph{difference differential equation},
$\ds \partial_z^k \varphi(z)=t(z)
+2e^{z/2}\varphi(z/2)$, whose solution involves
basic hypergeometric functions.
Catalan sums of the form $\sum_k \binom{2n}{n-k}\nu(k)$, with
$\nu$ an arithmetical function, arise in the statistics of ``order''
(also known as Horton-Strahler number) of trees
(Flajolet \& Prodinger,
1986). Both
cases are first subjected to a Mellin transform
analysis, which provides the relevant singular expansions.
Periodic fluctuations similar to the case of balanced trees then result
from singularity analysis.

A notable parallel to the  paradigm
of generating functions and singularity analysis has been developed by Vall\'ee
in a series of papers.
In her framework, singularities of certain transfer operators
(of Ruelle type)
replace singularities of generating functions.
See,~\emph{e.g.},~\cite{Vallee98b,Vallee01} for applications to
Euclidean algorithms and statistics on sequences produced by a general
model of dynamical
sources.

\section{Limit laws}

\vskip-5mm \hspace{5mm}

One of the important features of singularity analysis,
in contrast with Darboux's method or (real) Tauberian theory, is to allow
for \emph{uniform} estimates. This makes it possible to
analyse asymptotically
coefficients of multivariate generating functions, $f(z,u)$,
where the auxiliary
variable~$u$ marks some combinatorial parameter~$\chi$. One
\emph{first} proceeds
to extract $f_n(u):=[z^n]f(z,u)$  by considering $f(z,u)$ as a parameterized
family of univariate {\GF}s to which singularity analysis is applied.
(The coefficients $f_n(u)$ are, up to normalization,
probability generating functions of $\chi$.) A \emph{second} level of inversion
is then achieved by the standard continuity theorems for
probability characteristic functions (equivalently Fourier
transforms). Technically, consideration of a (small) neighbourhood of $u=1$ is
normally sufficient for extracting central limit laws.

Two important cases are those of a smoothly varying singularity and
of a smoothly varying exponent. In the first case, $f(z,u)$ has
a constant singular exponent~$\alpha_0$ and one has
$f(z,u)\sim c(u)(1-z/\rho(u))^{-\alpha_0}$.
Then, uniformity of singularity analysis implies
the estimate $f_n(u)/f_n(1)\sim (\rho(1)/\rho(u))^n$.
In other words, the probability generating function of~$\chi$ over
objects of size~$n$ is analytically similar to the {\GF} of a sum of
independent random variables---this situation is described as a
\emph{``quasi-powers'' approximation}. A Gaussian limit law for $\chi$ results
from the continuity theorem, with mean and variance that grow
in proportion to~$n$. The other case of a smoothly varying exponent
is dealt with similarly: one has
$f(z,u)\sim c(u)(1-z/\rho)^{-\alpha(u)}$ implying
 $f_n(u)/f_n(1)\sim n^{\alpha(u)-\alpha(1)}$; this is
once more a quasi-power
approximation, but with the parameter now in the scale
of $\log n$.  (See Gao \& Richmond, 1992, for hybrid cases.)

The technology above builds on early works of Bender~\cite{Bender73},
continued by Flajolet \& Soria~\cite{FlSo93,Soria90}, and
H.~K.~Hwang~\cite{Hwang94}.
In particular, under general conditions, the following
hold:
a local limit law expresses convergence to the Gaussian density;
speed of convergence estimates result from the Berry-Esseen inequalities;
large deviation estimates derive from singularity analysis applied at
fixed real values~$u\not=1$.

\begin{example}\label{poly-exa} \emph{Polynomials over finite fields}.
Consider the family $\cal P$ of all polynomials with coefficients in
the Galois field~$\mathbb{F}_q$.
A polynomial being determined by its sequence of coefficients,
the {\GF} $P(z)$ of all polynomials has a polar singularity.
Furthermore, the unique factorization property implies that $\cal P$ is
isomorphic to the class of all multisets ($\mathfrak{M}$) of the irreducible
polynomials~$\cal I$:  $\cal P\simeq \frak{M}\{\cal I\}$.
Since taking multisets corresponds to exponentiating singularities
of {\GF}s, the
singularity of the {\GF} $I(z)$ is logarithmic. By singularity analysis,
the number of irreducible polynomials is asymptotic to~$q^n/n$---this is
an analogue of the prime number theorem, which was already known to
Gau{\ss}. The bivariate {\GF} of the number of irreducible factors in polynomials
turns out to be
of the singular type $(1-qz)^{-u}$, with a smooth variable exponent,
so that:
\emph{the number of irreducible factors of a random polynomial over
$\mathbb{F}_q$ is asymptotically Gaussian}. This parallels the
Erd{\H o}s-Kac theorem for integers. Similar developments lead to a
complete analysis of  a major polynomial factorization algorithm
(Flajolet, Gourdon \& Panario, 2001).
\end{example}

Movable singularities and exponents occur frequently
in the analysis of parameters defined by recursion,
leading to algebraic or differential equations, which ``normally''
admit a smooth perturbative analysis.

\begin{example} \emph{Patterns in random strings.}
Let $\Omega$ be the total number of occurrences
of a fixed pattern (as a contiguous block)
 in a random string over a finite
alphabet. For either the Bernoulli model, where letters are
independently identically distributed, or the Markov model, the
bivariate {\GF}, with $z$ marking the length of the random string and $u$ the
number $\Omega$ of occurrences, is a \emph{rational function},
as it corresponds to a
finite-state device. Perron-Frobenius
properties apply for positive~$u$. Therefore the bivariate {\GF}
viewed as a function of~$z$ has a simple dominant pole at some
$\rho(u)$ that is an algebraic (and holomorphic) function of~$u$,
for $u>0$. Quasi-powers
approximations therefore hold and the limit law of $\Omega$
in random strings of length~$n$ is
Gaussian. Such facts holds for very general notions of patterns
and are developed systematically in Szpankowski's
book~\cite{Szpankowski01}.
\end{example}

\begin{example} \emph{Non-crossing graphs.}
Consider graphs with vertex set the $n$th roots of unity,
constrained to have only non-crossing edges; let the
parameter~$\chi$ be the number of connected components.
The bivariate {\GF} $G(z,u)$ is an \emph{algebraic function}
satisfying
\[G^3+(2w^3z^2-3w^2z+w-3)G^2+(3w^2z-2w+3)G+w-1=0
.
\]
%
$G(z,1)$ has a dominant singularity at $\rho(1)=3/2-\sqrt{2}$
which gets smoothly perturbed to $\rho(u)$ for $u$ near~1. The
singularity type is consistently of the form $(1-z/\rho(u))^{1/2}$.
A central limit law results for the number of
components in such
graphs~(Flajolet--Noy, 1999).
Drmota has given general conditions ensuring Gaussian laws for problems
similarly modelled by multivariate algebraic functions~\cite{Drmota97}.
\end{example}

\begin{example} \emph{Profile of quadtrees}.
Refer to Example~\ref{quad-exa}. The bivariate {\GF} $f(z,u)$
of node levels in quadtrees satisfies an equation, which, for
dimension~$d=3$ reads
\[
f(z,u)=1+2^3u\int_0^z \frac{dx_1}{x_1(1-x_1)}
\int_0^{x_1}\frac{dx_2}{x_2(1-x_2)}
\int_0^{x_2} f(x_3,u)\frac{dx_3}{1-x_3}.
\]
This corresponds to a \emph{linear differential equation} with coefficients
in $\mathbb{C}(z,u)$ and a fixed singularity at~$z=1$.
The indicial equation is an algebraic one parameterized by~$u$
and, when
$u\approx1$, there is a unique largest branch $\alpha(u)$ that determines
the dominant regime of the form~$(1-z)^{-\alpha(u)}$. This is the case of
a movable exponent inducing a central limit law: \emph{The level profile of
a $d$-dimensional quadtree is asymptotically Gaussian.}
Such properties are expected in general for models that are
perturbations of linear differential equations with a fixed Fuchsian
singularity (Flajolet \& Lafforgue, 1994).
\end{example}

Finally,
singularity analysis also intervenes by making it possible to
``pump'' moments of combinatorial distributions.
Examples include the height of
trees discussed earlier, as well as tree path length
(Louchard 1983, Tak{\'a}cs 1991) and the
construction cost of hashing tables (Flajolet, Poblete \& Viola,
1998). The latter problems were first shown in this way to
converge to Brownian Excursion Area.

\section{Conclusions}

\vskip-5mm \hspace{5mm}

Elementary combinatorial structures are enumerated by generating
functions that satisfy a rich variety of functional relations.
However, the singular types that are observed are usually somewhat
restricted, and \emph{driven by combinatorics}.
In simpler cases, the generating functions are explicit
combinations of a standard set of special functions.
Next, implicitly defined functions (associated with recursion)
have singularities that
arise from failures of the implicit function theorem
and are consequently of the algebraic type, often
with exponent $\frac12$.
Linear differential equations have a well-established
classification theory that, in the Fuchsian case, leads to
algebraic-logarithmic
singularities. In all such cases, the singular expansion is known to
be valid outside of the original disc of convergence of the generating
function. This means that singularity analysis is \emph{automatically}
applicable, and precise asymptotic expansions of coefficients result.

Parameters of combinatorial structures, provided they remain
``simple'' enough, lead to local deformations (via an auxiliary
variable~$u$ considered near~1)
of the functional relations
defining univariate counting generating functions. Under fairly general
conditions, such deformations are amenable to perturbation theory
and admit of uniform expansions near
singularities. Then, since the singularity analysis process preserves
uniformity, limit laws result via the continuity theorem for
characteristic functions. In this way, the behaviour of a large number of
parameters of elementary combinatorial structures becomes
predictable. (The theory of functions of several complex
variables is thus bypassed. See Pemantle's recent
work~\cite{Pemantle00} based on this theory
for a global characterization of all the asymptotic regimes involved.)

The generality of the singular approach makes it
even possible to discuss \emph{combinatorial schemas}
at a fair level of generality~\cite{FlSe01bsp,FlSo93,Hwang94,Soria90},
the case of polynomial factorization (Ex.~\ref{poly-exa}) being typical.
Roughly, combinatorial constructions viewed as ``singularity
transformers'' dictate asymptotic regimes and probabilistic laws.
Analytic combinatorics then represents an attractive alternative to
probabilistic methods, whenever a strong analytic
structure is present---this is the case
for most combinatorial problems that are ``decomposable''
and amenable to the generating function methodology.
Very precise asymptotic
information on the randomness properties of large random
combinatorial objects results from there. This in turn has useful
implications in the analysis of many fundamental algorithms
and data structures of computer science, following the steps of
Knuth's pioneering works~\cite{Knuth00}.

\providecommand{\bysame}{\leavevmode\hbox to3em{\hrulefill}\thinspace}

\label{lastpage}


\begin{thebibliography}{10}

\bibitem{Bender73}
Edward~A. Bender, \emph{Central and local limit theorems applied to asymptotic
  enumeration}, Journal of Combinatorial Theory \textbf{15} (1973), 91--111.

\bibitem{Bender74}
\bysame, \emph{Asymptotic methods in enumeration}, SIAM Review \textbf{16}
  (1974), no.~4, 485--515.

\bibitem{Drmota97}
Michael Drmota, \emph{Systems of functional equations}, Random Structures \&
  Algorithms \textbf{10} (1997), no.~1--2, 103--124.

\bibitem{FiFlKa02}
James~A. Fill, Philippe Flajolet, and Nevin Kapur, \emph{{Singularity Analysis,
  Hadamard Products, and Trees Recurrences}}, Preprint, 2002.

\bibitem{Flajolet99}
Philippe Flajolet, \emph{Singularity analysis and asymptotics of {Bernoulli}
  sums}, Theoretical Computer Science \textbf{215} (1999), no.~1-2, 371--381.

\bibitem{FlOd82}
Philippe Flajolet and Andrew~M. Odlyzko, \emph{The average height of binary
  trees and other simple trees}, Journal of Computer and System Sciences
  \textbf{25} (1982), 171--213.

\bibitem{FlOd90b}
\bysame, \emph{Singularity analysis of generating functions}, SIAM Journal on
  Algebraic and Discrete Methods \textbf{3} (1990), no.~2, 216--240.

\bibitem{FlSe01bsp}
Philippe Flajolet and Robert Sedgewick, \emph{Analytic combinatorics}, 2001,
  Book in preparation; see also INRIA Research Reports 1888, 2026, 2376, 2956,
  3162, 4103.

\bibitem{FlSo93}
Philippe Flajolet and Mich{\`e}le Soria, \emph{General combinatorial schemas:
  {G}aussian limit distributions and exponential tails}, Discrete Mathematics
  \textbf{114} (1993), 159--180.

\bibitem{Ford60sp}
W.~B. Ford, \emph{Studies on divergent series and summability}, 3rd ed.,
  Chelsea, New York, 1960, (From two books originally published in 1916 and
  1936.).

\bibitem{Hwang94}
Hsien-Kuei Hwang, \emph{Th\'eor\`emes limites pour les structures combinatoires
  et les fonctions arithmetiques}, Ph.D. thesis, {\'E}cole Polytechnique,
  December 1994.

\bibitem{Knuth00}
Donald~E. Knuth, \emph{Selected papers on analysis of algorithms}, CSLI
  Publications, Stanford, CA, 2000.

\bibitem{Mahmoud92}
Hosam~M. Mahmoud, \emph{Evolution of random search trees}, John Wiley, New
  York, 1992.

\bibitem{Odlyzko82}
A.~M. Odlyzko, \emph{Periodic oscillations of coefficients of power series that
  satisfy functional equations}, Advances in Mathematics \textbf{44} (1982),
  180--205.

\bibitem{Odlyzko95}
\bysame, \emph{Asymptotic enumeration methods}, Handbook of Combinatorics
  (R.~Graham, M.~Gr{\"o}tschel, and L.~Lov{\'a}sz, eds.), vol.~{II}, Elsevier,
  Amsterdam, 1995, 1063--1229.

\bibitem{Pemantle00}
Robin Pemantle, \emph{Generating functions with high-order poles are nearly
  polynomial}, Mathematics and computer science (Versailles, 2000),
  Birkh\"auser, Basel, 2000, 305--321.

\bibitem{Polya37}
G.~P{\'o}lya, \emph{{K}ombinatorische {A}nzahlbestimmungen f{\"u}r {G}ruppen,
  {G}raphen und chemische {V}erbindungen}, Acta Mathematica \textbf{68} (1937),
  145--254.

\bibitem{Salvy91}
Bruno Salvy, \emph{Asymptotique automatique et fonctions g{\'e}n{\'e}ratrices},
  Ph. {D}. thesis, {\'E}cole Polytechnique, 1991.

\bibitem{Soria90}
Mich{\`e}le Soria-Cousineau, \emph{M{\'e}thodes d'analyse pour les
  constructions combinatoires et les algorithmes}, Doctorat {\`es} sciences,
  Universit\'e de Paris--Sud, Orsay, July 1990.

\bibitem{Szpankowski01}
Wojciech Szpankowski, \emph{Average-case analysis of algorithms on sequences},
  John Wiley, New York, 2001.

\bibitem{Vallee98b}
Brigitte Vall{\'e}e, \emph{Dynamics of the binary {E}uclidean algorithm:
  Functional analysis and operators}, Algorithmica \textbf{22} (1998), no.~4,
  660--685.

\bibitem{Vallee01}
\bysame, \emph{Dynamical sources in information theory: Fundamental intervals
  and word prefixes}, Algorithmica \textbf{29} (2001), no.~1/2, 262--306.

\bibitem{WiZe85b}
Jet Wimp and Doron Zeilberger, \emph{Resurrecting the asymptotics of linear
  recurrences}, Journal of Mathematical Analysis and Applications \textbf{111}
  (1985), 162--176.

\bibitem{Zeilberger90}
Doron Zeilberger, \emph{A holonomic approach to special functions identities},
  Journal of Computational and Applied Mathematics \textbf{32} (1990),
  321--368.

\end{thebibliography}
\end{document}